\documentclass[11]{siamltex}
\usepackage{amsmath,amssymb,amsfonts,graphicx,fancybox,shadow,color}
\usepackage{hyperref}

\newcommand{\RR}{{\mathbb{R}}}
\newcommand{\CC}{{\mathbb{C}}}

\def\tu {\tilde{u}}
\def\tF {\tilde{F}}

\def\uu{{\mathbf u}}
\def\UU{{\mathbf U}}

\def\be#1\ee{\begin{equation}#1\end{equation}}

\begin{document}
\bibliographystyle{plain}

\pagestyle{myheadings}

\title{On extension of the data driven  ROM  inverse scattering framework to partially nonreciprocal arrays }
\author{
 V. Druskin\footnotemark[1],  S. Moskow\footnotemark[2] and M. Zaslavsky\footnotemark[3]}

\renewcommand{\thefootnote}{\fnsymbol{footnote}}

\footnotetext[1]{Worcester Polytechnic Institute, Department of Mathematical Sciences,
Stratton Hall,
100 Institute Road, Worcester MA, 01609 (vdruskin1@gmail.com)}
\footnotetext[2]{Department of Mathematics, Drexel University, Korman Center, 3141 Chestnut Street, Philadelphia, PA 19104
(moskow@math.drexel.edu)}
\footnotetext[3]{Schlumberger-Doll Research Center, 1 Hampshire St., 
Cambridge, MA 02139-1578 (mzaslavsky@slb.com)}

\maketitle

\begin{abstract} \ Data-driven reduced order models (ROMs) recently emerged as powerful tool for the solution of  inverse scattering problems. The main drawback of this approach is that it was limited to the measurement arrays with reciprocally collocated  transmitters and receivers, that is, square symmetric matrix (data) transfer functions. To relax this limitation,  we  use  our   previous work \cite{DrMoZa}, where the ROMs  were combined with the  
Lippmann-Schwinger integral equation to produce a direct nonlinear inversion method.  In this work we extend this approach to more general  transfer functions, including those that are non-symmetric, e.g., obtained by adding only  receivers or sources.    The ROM is constructed based on the symmetric subset of the data and is used to construct all internal solutions. Remaining receivers are then used directly in the Lippmann-Schwinger equation. We demonstrate the new approach on a number of 1D and 2D examples with non-reciprocal arrays, including   a single input/multiple outputs (SIMO) inverse problem, where the data is given by just a single-row matrix transfer function.\end{abstract}

\section{Introduction}

In this work we extend the reduced order model (ROM) approach which was used previously for inverse impedance, scattering and diffusion \cite{borcea2011resistor,borcea2014model,druskin2016direct,druskin2018nonlinear,borcea2017untangling,borcea2019robust,BoDrMaMoZa,borcea2020reduced, DrMoZa,  Borcea2021ReducedOM} to nonsymmetric data sets. Here, we consider the inverse diffusion problem given spectral data, although the technique to handle non-symmetric transfer functions can be applied to a much wider class of problems. 

In the ROM framework for solving multidimensional inverse problems, the ROM is constructed from a symmetric data set, that is, for coinciding source and receiver pairs. That is, the ROM is chosen precisely to match this data set, see
\cite{druskin2016direct,druskin2018nonlinear,borcea2017untangling,borcea2019robust, BoDrMaMoZa, Borcea2021ReducedOM, borcea2021reduced}.  Then, the ROM is transformed to a sparse form (tridiagonal for single input single output (SISO) problems, block tridiagonal for multiple input/output (MIMO) problems) by Lanczos orthogonalization. 
The data-driven ROM, in this orthogonalized form, can be viewed as a discrete network which has entries for which their dependence on the unknown pde coefficients is approximately linear \cite{borcea2011resistor,druskin2016direct,borcea2014model,borcea2020reduced,Borcea2021ReducedOM}.  
That is, the main nonlinearity of the inverse problem was absorbed by the orthogonalization process. 
This process is related to the seminal works of Marchenko, Gelfand, Levitan and Krein on inverse spectral problems, and to the idea of spectrally matched second order staggered finite-difference grids first introduced in \cite{druskin1999gaussian} and first used for direct inversion in \cite{BoDr}.  
%

The data-driven ROM can be viewed as a projection operator, or Galerkin system \cite{doi:10.1137/1.9781611974829.ch7}. 
A crucial property, which was first noticed in \cite{druskin2016direct}, is that the Galerkin basis corresponding to the unknown coefficient is very close to the one from the homogeneous problem. This led to the back projection algorithm \cite{druskin2018nonlinear}, which allowed for direct reconstructions in multiple dimensions. In \cite{BoDrMaMoZa}, it was found that thanks to this weak dependence of the basis functions on the unknown coefficients, the ROM can also be directly used to generate internal solutions from boundary data only.

In the work \cite{DrMoZa}, the data-generated internal solutions ${\bf u}_p$ (corresponding to unknown coefficient $p$) introduced in \cite{BoDrMaMoZa} were used in the Lippmann-Schwinger integral equation framework. If we consider the data $F_p$ for the unknown coefficient $p$ and background data $F_0$ corresponding to a known background, then 
the Lippmann-Schwinger integral equation with respect to the unknown $p$ can be written as \be\label{eq:LipSwi} F_p -F_0 =-\langle u_0,p u_p\rangle \ee
where $\langle ,\rangle$  is the continuous $L^2$ inner product on the PDE domain, and where $u_p$ and $u_0$ are the unknown and background internal solutions respectively. For the Lippmann-Schwinger Lanczos integral equation {(LSL IE)}, we use the data generated internal solution ${\bf u}_p$ in place of $u_p$:
\be\label{eq:LipSwiL} F_p -F_0  \approx -\langle u_0,p {\bf u}_p\rangle.\ee
Since ${\bf u}_p$ is precomputed directly from the data without knowing $p$, (\ref{eq:LipSwiL}) becomes linear with respect to $p$. 

The Lippmann-Schwinger-Lanczos approach resolves directly one of the main limitations of the above mentioned works;  that it applies only to symmetric data sets. This is the subject of this manuscript.   Indeed, while the ROM must still be constructed from a symmetric subset of the data, the Lippmann-Schwinger Lanczos equation (\ref{eq:LipSwiL}) can be directly applied to any additional receivers, using the same internal solutions. Furthermore, additional sources (with receivers in the symmetric part) can also be added thanks to reciprocity.  Adding this additional data increases the range of the   Lippmann-Schwinger-Lanczos linear operator on the r.h.s. of \eqref{eq:LipSwiL}, thus  improving the quality of the solution of the inverse problem.

The structure of our extended multi-input/multi-output (MIMO)  measurement array can be summarized with the help of \eqref{eq:ext.matrix}.  The column and row numbers correspond respectively to the indices of the receivers (outputs) and transmitters (inputs). The conventional data-driven ROM requires collocated receiver and transmitter arrays (the symmetric data set), so the measurements in this case are given by a square matrix, that is,  the left upper block of  \eqref{eq:ext.matrix}.  We assume reciprocity, implying symmetric transfer functions, thus it is sufficient just to measure their upper triangular parts.  The extended receivers are shown in the upper right block. A common case is when all added receivers measure responses from one of the transmitters, resulting in adding a column.  However,  added receivers can potentially just measure responses from certain transmitters,  resulting in an arbitrary sparsity pattern in the upper right block.  Likewise, one can add transmitters by adding elements to the  lower left  block. By reciprocity, they can be symmetrically reflected into the right upper one. 

\be\label{eq:ext.matrix} 
	\begin{pmatrix} \begin{pmatrix}.& . & . & . & . \\
    \  & . & . & . & . \\
     \  &  \  & . & . & . \\
     \  &  \  &  \  & . & . \\
     \  &  \  &  \  &  \  & . \\ 
	\end{pmatrix}  & \begin{matrix}.& \ & \ & \ & . \\
    .  & \ &\  & \  & \  \\
     .  &  \  &  \  &  .  &  \  \\
     .  &  \  &  .  & \  &  \  \\
     . &  \  &  \  &  .  &  \  \\ 
	\end{matrix}   \\  
	 &  &   \\ & & \\  &   &  
	\\   & &  
	\end{pmatrix}
\ee

This paper is organized as follows. In Section 2 we describe the entire process in detail for a one dimensional single input multiple output (SIMO) problem. We briefly describe the construction of the ROM from the single input/  single output (SISO)  part of data, the Lanczos orthogonalization process and the generation of the internal solutions.  We then show its use in the Lippmann-Schwinger equation for the full SIMO data. 
The generalization of this process to SIMO or nonsymmetric  MIMO arrays in higher dimensions is described in Section 3. Section 4 contains numerical experiments.

\section{One dimensional SIMO problem}
We begin this work with a one dimensional problem to demonstrate the approach.  In the first subsection we describe the problem setup one source and two receivers, one of which is at the source.  We note the transfer function for this problem has the SISO transfer function as a component, which we use to construct the ROM. We review this briefly in the second subsection, along with tridiagonalization  of the ROM and generation of the internal solutions.  In the last subsection we show how to use the full transfer function and the internal solutions in the Lippmann-Schwinger equation in order to solve the fully nonlinear SIMO inverse problem. 

  \subsection{ Description of the SIMO problem} We start by considering the following single input multiple output (SIMO) inverse problem in one dimension
	\be\label{eq:1D} -\frac{d^2u}{dx^2}(x,\lambda) +p(x)u(x,\lambda)+\lambda u(x,\lambda)=g^{(1)}(x) \quad \frac{du}{dx}|_{x=0}=0,\ \frac{du}{dx}|_{x=L}=0,
	\ee
where $0<L\le\infty$. The source $g^{(1)}(x)$ is assumed to be a compactly supported real distribution localized near the origin, for example, roughly speaking, $ g^{(1)} = \delta(x-\epsilon)$ with small $\epsilon>0$.  We also use a similar distribution $g^{(2)}$ near the right boundary $x=L$ to represent the second receiver. 
Consider  $\lambda_j\in \CC\setminus \RR_-$, $j=1,\ldots, m$, since for nonnegative $p$ the above resolvent is well defined for $\lambda$ off of the negative real axis.
We note that the above formulation can be obtained via the Laplace transform of the one dimensional diffusion or wave problem.
The SIMO transfer function is then 
\begin{eqnarray} \label{eq:transfer}
	F(\lambda) &=& \left[ \int_0^L g^{(1)}(x) u(x,\lambda)dx , \int_0^L g^{(2)}(x) u(x,\lambda)dx \right] \\ &=& 
	\left[ \langle g^{(1)},u\rangle ,  \langle g^{(2)},u\rangle\right] \end{eqnarray} 
where throughout the paper we use $\langle, \rangle$ to denote the continuous Hermitian $L^2$ inner product $$\langle w,v\rangle=\int_0^L \bar{w}(x)v(x)dx ,$$ which in the one dimensional case is $L^2(0,L)$.

For simplicity of exposition we consider  $2m$ real spectral data points that is, we consider the data 
\begin{equation} \label{realdata} F(\lambda)|_{\lambda=\lambda_j}\in\RR, \ \ \ \frac{dF(\lambda)}{d\lambda}|_{\lambda=\lambda_j}\in \RR \ \ \ \mbox{for} \ \  j=1,\ldots, m.\end{equation}  
In this case all solutions are real and the conjugates are unnecessary. 
For complex data points, all of the following holds, see \cite{DrMoZa}.   The SIMO inverse problem is then to determine $p(x)$ in (\ref{eq:1D}) from the data  (\ref{realdata}).

\subsection{Construction of the data-driven ROM, orthogonalization and internal solutions}

The first component of the transfer function (\ref{eq:transfer}) is a SISO transfer function, 
$$F^\sigma := ( F(\lambda) )_1 $$
is the symmetric portion of $F$, and can be used to construct the ROM exactly as in previous works e.g. \cite{DrMoZa}.  We describe this now briefly. The exact solutions to (\ref{eq:1D}) , $\{ u(x,\lambda_j)\}$, form a basis for the projection subspace $$\mathbb{U}=\mbox{span}\{u_1(x)=u(x,\lambda_1),\ldots , u_m(x)=u(x,\lambda_m)\}.$$  We define the data-driven ROM as the Galerkin system corresponding to $\mathbb{U}$ \be\label{eq:ROM}
(S+\lambda M)c=b
\ee
where $S,M\in\mathbb{R}^{m\times m}$ are symmetric positive definite matrices with the stiffness matrix $S$ given by $$S_{ij}=\langle {u}'_i,u'_j\rangle +\langle p{u}_i,u_j\rangle $$ and mass matrix $M$ given by  $$M_{ij}=\langle {u}_i,u_j\rangle. $$ The right hand side $b\in\mathbb{R}^m$ is a column vector with components $$b_j =\langle  {u}_j, g  \rangle, $$ and the Galerkin solution for the system is determined by the  vector   $c\in\mathbb{R}^{m}$ depending on $\lambda$.  Note that $c= c(\lambda)$ corresponds to a column vector of coefficients of the solution with respect to the above basis of exact solutions. The matrices $S$ and $M$ can be obtained directly from the data, without knowing the exact solutions, from the formulas
\begin{equation}
\label{eq:massmtr}
M_{ij}=\frac{{F^\sigma}(\lambda_i)-F^\sigma(\lambda_j)}{\lambda_j-{\lambda}_i}, \ \ \ M_{ii} = -{dF^\sigma\over{d\lambda}}(\lambda_i).
\end{equation}
and
\begin{equation}
\label{eq:stifmtr}
S_{ij}=\frac{F^\sigma(\lambda_j)\lambda_j-{F^\sigma}(\lambda_i){\lambda}_i}{\lambda_j-{\lambda}_i}, \ \ \ S_{ii} = {d(\lambda F^\sigma)\over{d\lambda}}(\lambda_i).
\end{equation}
Due to the matching conditions,  the ROM transfer function corresponding to this Galerkin system   
     \[ \tF^\sigma(\lambda)=b^\top c\] 
matches (the symmetric part of) the data exactly. This is well known; a proof of this for real $\lambda_i$ is in \cite{BoDrMaMoZa}, and for complex $\lambda$ in \cite{DrMoZa}. 
Furthermore, for any $\lambda$, the solution to (\ref{eq:1D}) is close to its Galerkin projection 
\[u(\lambda)\approx \tu(\lambda)={V}c ={V} (S+\lambda M)^{-1}b \]
where  ${V}$  represents the row vector of basis functions $u_i$,
$${ V} = [ u_1 ,\ldots , u_m ] .$$

Next we perform a change of basis, that is, we orthogonalize by using the Lanczos algorithm. More precisely, we run $m$ steps of the { $M$-symmetric}  Lanczos algorithm corresponding to  matrix $A=M^{-1}S$ and initial vector  $M^{-1}b$. This yields tridiagonal matrix $T\in\RR^{m\times m}$ and $M$-orthonormal Lanczos vectors $q_i\in \CC^m$, such that
	\be\label{eq:lancz} AQ =Q T , \qquad Q^\top MQ=I,\ee
where $$Q=[q_1, q_2, \ldots, q_m]\in{\CC^{m\times m}},$$ and $$q_1=M^{-1}b/\sqrt{b^\top M^{-1}b}.$$
 The new basis is orthonormal in $L^2(0,L)$ and is given by the vector $${VQ=  [ \sum_{j=1}^mq_{j1}u_j ,\ldots , \sum_{j=1}^m q_{jm}u_j ]}  .$$   The Galerkin solutions and transfer function can then be written in this new basis as
	\be \label{eq:state} \tu(\lambda)=\sqrt{b^\top M^{-1}b}V Q(T+\lambda I)^{-1}e_1, \ee \be \tF^\sigma(\lambda)= (b^\top M^{-1}b)  e_1^\top (T+\lambda I)^{-1}e_1\ee
where $e_1 = (1,0,\ldots,0)^T $ is the first coordinate column vector in $\RR^m$.

Now, we use this orthogonalized ROM to produce internal solutions. Of course, we do not know $p$, so we don't know the original basis of exact solutions $V$. What we do instead is to replace the unknown orthogonalized internal solutions $VQ$ with orthogonalized background solutions $V_0Q_0$ corresponding to background $p_0=0$. Here $V_0$ is the row vector of background solutions 
$$V_0 = [ u^0_1,\ldots , u^0_m ] $$ 
to (\ref{eq:1D}) corresponding to $p=p_0=0$ and the same spectral points 
$\lambda= \lambda_1, \ldots \lambda_m$. A ROM for this background problem is computed in the same way, and $Q_0$ is computed from its Lanczos orthogonalization.  That is, one can compute an approximation $\uu$ to the unknown internal solution 
$u(x,\lambda)$ using
\begin{equation}\label{internal1d}  u \approx \uu= \sqrt{b^\top M^{-1}b} V_0 Q_0(T+\lambda I)^{-1}e_1  \end{equation}
which is obtained from data only.

\subsection{Nonlinear inverse problem} 
We now use the Lippmann-Schwinger formulation to solve the nonlinear inverse problem for all of the data, not just the symmetric part, which is the SISO transfer function in this case.  
From  (\ref{eq:transfer}) and its background counterpart we obtain the Lippmann-Schwinger equation 
	\be  \label{eq:int}F_0(\lambda_j)-F(\lambda_j)=\left[ \int {u}^{(1)}_0(x,\lambda_j ) u(x,\lambda_j )p(x)dx, \int {u}^{(2)}_0(x,\lambda_j ) u(x,\lambda_j )p(x)dx\right] \ee
for $j=1,\ldots, m$
where $u_0^{(i)}$ is the solution to the background problem
\be\label{eq:1D0} -\frac{d^2u_0^{(i)}}{dx^2}(x,\lambda) +\lambda u_0^{(i)}(x,\lambda)=g^{(i)}(x)
	\ee with zero Neumann conditions. 
Correspondingly, $F_0$ is the background transfer function
\begin{equation} \label{eq:backtransfer}
F_0	=\left[ \langle g^{(1)},u^{(1)}_0\rangle ,  \langle g^{(2)},u^{(1)}_0\rangle\right]. \end{equation} 
For real $\lambda_j\in \RR$ we then have $2\times2m$ equations 
\begin{multline}\label{eq:intd}F_0(\lambda_j) -F(\lambda_j) \\ =\left[ \int  u^{(1)}_0(x,\lambda_j) u(x,\lambda_j )p(x)dx,\int  u^{(2)}_0(x,\lambda_j) u(x,\lambda_j )p(x)dx\right] ,\end{multline} \begin{multline}   \frac{d}{d\lambda}(F_0-F)|_{\lambda= \lambda_j} \\ =\left[\int \frac{d}{d\lambda} [u^{(1)}_0(x,\lambda) u(x,\lambda)]_{\lambda=\lambda_j} p(x)d  x, \int \frac{d}{d\lambda}u^{(2)}_0(x,\lambda) u(x,\lambda)]_{\lambda=\lambda_j} p(x)d  x\right] \label{eq:intdd} \end{multline}
for $j=1,\ldots,m$. 
If we put the background solutions in the vector
\begin{equation} { U}_0 = \left[ u^{(1)}_0 , u^{(2)}_0\right] \nonumber\end{equation}
we can write (\ref{eq:intd}) and (\ref{eq:intdd}) as 
\begin{equation}\label{eq:intd2}F_0(\lambda_j) -F(\lambda_j) = \int {  U}_0 u(x,\lambda_j )p(x)dx,  \end{equation} \begin{equation}   \frac{d}{d\lambda}(F_0-F)|_{\lambda= \lambda_j} =\int \frac{d}{d\lambda} \left[ {U}_0(x,\lambda) u(x,\lambda)\right]_{\lambda=\lambda_j} p(x)d  x \label{eq:intdd2} \end{equation}
for $j=1,\ldots,m$. 
As usual, the internal solutions $u(x,\lambda_j )$ and their derivatives with respect to $\lambda$ are unknown, and depend on $p$, so the system (\ref{eq:int}-\ref{eq:intd}) is nonlinear with respect to $p$. Now, just as in \cite{DrMoZa} we replace $u(x,\lambda )$  in (\ref{eq:int}-\ref{eq:intd}) with its approximation $$ u \approx \uu=\sqrt{b^*M^{-1}b}V_0 Q_0 (T+\lambda I)^{-1}e_1,$$ 
which we obtained from the ROM constructed from the symmetric part of the data.
We can write the new system for $p$ as
\be\label{eq:oper}
\delta F= \int W(x) p(x) dx
\ee
where $$\delta F=[(F_0-F)(\lambda_1),\ldots, (F_0-F)(\lambda_m), {d\over{d\lambda}}(F_0-F)(\lambda_1),\ldots, {d\over{d\lambda}}(F_0-F)(\lambda_m)]\in\mathbb{R}^{4m},$$  and     $$W=[\uu {U}_0(\lambda_1),\ldots,
\uu {U}_0(\lambda_m),{d\over{d\lambda}}(\uu{U}_0)|_{\lambda=\lambda_1},\ldots,
{d\over{d\lambda}}(\uu {U}_0)|_{\lambda=\lambda_m}]$$ 
is a $4m$-dimensional vector of functions on $(0,L)$. Recall $\uu$  is computed directly from the symmetric part of the data without knowing $p $, by using (\ref{internal1d}),  making nonlinear system (\ref{eq:oper})  linear. We  continue to refer to (\ref{eq:oper}) as a Lippmann-Schwinger-Lanczos system. 
{ Setup \eqref{eq:oper}  corresponds to the full $1\times 2$ array in \eqref{eq:ext.matrix}.}

\section{Multidimensional MIMO problem}
In this section we consider a more general multidimensional case where we may have any number of sources and receivers. The symmetric part of the transfer function will consist of all of the data from coinciding source/receiver pairs,  that is, it corresponds to the  left upper block of
\eqref{eq:ext.matrix}. Again, the symmetric part will be used to build the ROM as in \cite{DrMoZa}, while the remaining data will be used later in the Lippmann-Schwinger formulation. We will assume that the nonsymmetric, or remaining part of the transfer function consists of receivers, { that is, they will be located in the right upper block of \eqref{eq:ext.matrix}.  As it was already mentioned in the introduction, } this is without loss of generality, since any remaining sources can equivalently be viewed as receivers by reciprocity. 
\subsection{Description of the MIMO problem}
{We consider a formulation that can be obtained via the Laplace transform of the diffusion or wave problem, which can be written as} the following boundary value problem on $\Omega\in\RR^d$:
\be\label{d-Schrod}-\Delta u^{(r)} +p u^{(r)} +\lambda u=g^{(r)}, \quad  \frac{d u}{d\nu}\large|_{\partial \Omega}=0, \  r=1,\ldots, K,\ee
where $g^{(r)}$ are $K$ localized  sources, e.g., boundary charge distributions,  supported near or at an accessible part $S$ of $\partial \Omega$. Consider also a set of distributions representing our receivers, which will include all of the sources distributions. That is, let
$$G=[g^{(1)}, g^{(2)}, \ldots, g^{(K)}, \ldots, g^{(L)}],$$ where $L\ge K$. Let $$U=[u^{(1)}, u^{(2)}, \ldots, u^{(K)}]$$
be our solution functions corresponding to the $K$ sources. Then the multiple-input multiple output (MIMO) transfer function is a $K\times L$ matrix valued function of $\lambda$
	\be\label{eq:transferMIMO}
	F(\lambda)= \langle G,U\rangle \in\mathbb{R}^{K\times L}
	\ee
where $\langle, \rangle$ represents the continuous $L^2(\Omega)$ inner product, that is, 
$$\langle G,U\rangle = \int_\Omega U^\top G dx, $$ is matrix valued with 
$$(F(\lambda))_{ij} = \int_\Omega u^{(i)} g^{(j)}dx. $$  We consider the inverse problem with data given by $2m$ real symmetric $K\times L$ matrices, that is, $$F(\lambda)|_{\lambda=\lambda_j}\in\RR^{K\times L},$$ and  $$\frac{F(\lambda)}{d\lambda}|_{\lambda=\lambda_j}\in \RR^{K\times L},$$
for real spectral points $\{ \lambda_j\}_{j=1,\ldots m}$, while we note again that complex data will follow in the same way \cite{DrMoZa}.

\subsection{Construction of the data-driven ROM, block Lanczos, and internal solutions}
We now construct the ROM based on the symmetric part of the transfer function:
$$ F^\sigma (\lambda) := (F(\lambda))|_{i,j=1\ldots K}$$
and correspondingly define $$G^\sigma = [g^{(1)}, g^{(2)}, \ldots, g^{(K)}].$$
We consider the $mK$ dimensional projection subspace $$\mathbb{U}=\mbox{span}\{U_1(x)=U(x,\lambda_1),\ldots , U_m(x)=U(x,\lambda_m)\}$$ and define the MIMO data-driven ROM as the system \be\label{eq:ROMMIMO}
(S+\lambda M)C  = B
\ee
where $S,M\in\RR^{mK\times mK}$ are symmetric positive definite matrices,  $B\in\RR^{mK \times K}$ is the symmetric part of the data $$ B= \langle G^\sigma , U \rangle, $$ and the system solution $C=C(\lambda) \in\RR^{mK\times K} $ is a matrix valued function of $\lambda$, again corresponding to coefficients of the solution with respect to the above basis of exact solutions. Stiffness and mass matrices are block versions of (\ref{eq:ROM}), with blocks given by  $$S=(S_{ij}=\langle {U'}_i,U'_j\rangle)+\langle p{U}_i,U_j\rangle)$$ and $$M=(M_{ij}=\langle {U}_i,U_j\rangle).$$  Again $S$ and $M$ are obtained by imposing the conditions that the ROM transfer function  $\tF^\sigma(\lambda)=B^\top C $
matches the data.  Next, we perform block Lanczos tridiagonalization, with matrix $$A=M^{-1}S$$ and initial block vector $M^{-1}B$.  From this we obtain the block-tridiagonal matrix $$T\in\RR^{Km\times Km}$$ with  $K\times K$ blocks,  and  block-vectors $q_i\in \RR^{mK\times K}$ which are orthonormal with respect to the $M$ inner product, which again corresponds to continuous $L^2$ orthonormality. From this we obtain the block counterpart of (\ref{eq:lancz}) 
	$$ Q=[q_1, q_2,\ldots, q_m]\in{\RR^{mK\times mK}},$$ where $$ q_1=M^{-1}B(B^\top M^{-1}B)^{-1/2}.$$
The Galerkin projection of the true solution or state solution can then be written in the new Lanczos basis as
	\be \label{eq:stateB} \tilde{U}(\lambda)=\sqrt{B^\top M^{-1}B}V Q(T+\lambda I)^{-1}E_1,\ee
where $E_1\in\RR^{mK\times K}$ consists of the first $K$ columns of identity matrix $I\in\RR^{mK\times mK}$.
 We express the orthogonalized basis as 
$$VQ=  [ \sum_{j=1}^m U_jq_{j1} ,\ldots , \sum_{j=1}^m U_jq_{jm}] $$
where $\left\{q_{ji}\in \RR^{K\times K}\right\}^{m}_{i,j=1}$ are the blocks of matrix $Q$ and $U_j=[u^{(1)}_{j}(x),\ldots u^{(K)}_{j}(x)]$ are the corresponding vectors of solutions. Note that this basis depends on the unknown internal solutions, so as before we replace the basis $VQ$ in (\ref{eq:stateB}) with its background counterpart $V_0Q_0$. 
 This yields the data generated vector of internal solutions
	\be\label{eq:internMIMO} \UU := \sqrt{B^\top M^{-1}B}V_0 Q_0(T+\lambda I)^{-1}E_1.
	\ee
	
\subsection{Nonlinear MIMO inverse problem}
We now solve the fully nonlinear inverse problem using the entire non-symmetric data set.
From  (\ref{eq:transferMIMO}) and the equation for the background solutions we obtain the Lippman-Schwinger formulation
 \be\label{eq:intMIMO}F_0(\lambda_j)-F(\lambda_j)=\int U^\top (x,\lambda_j )U_0(x,\lambda_j ) p(x)dx, \qquad j=1,\ldots, m.\ee
Here $U_0$ represents all $L$ background solutions with $p=0$.
For $\lambda_j\in \RR$ we will have
\begin{eqnarray} \label{eq:intdMIMO}(F_0-F)|_{\lambda_j} &=& \int {U}^\top(x,\lambda_j )U_0(x,\lambda_j) p(x) dx, \nonumber \\  \frac{d}{d\lambda}(F_0-F)|_{\lambda_j}  &=& \int \frac{d}{d\lambda}( {U}^\top(x,\lambda_j )U_0(x,\lambda) )|_{\lambda=\lambda_j} p(x)d x.\end{eqnarray}
Now we replace $U(x,\lambda )$ with its approximation $\UU$ in (\ref{eq:intMIMO}) and (\ref{eq:intdMIMO}). 
Again, precomputing $\UU$ from the data only (\ref{eq:internMIMO}) will yield the linear system for $p$, {which is given by the upper triangular part of \eqref{eq:intdMIMO2} and \eqref{eq:intdMIMO2d}}:
\begin{eqnarray} \label{eq:intdMIMO2}(F_0-F)|_{\lambda_j} &=& \int \mathbf{U}^\top(x,\lambda_j )U_0(x,\lambda_j ) p(x) dx,  \\ {\label{eq:intdMIMO2d} }\frac{d}{d\lambda}(F_0-F)|_{\lambda_j}  &=& \int \frac{d}{d\lambda}(\mathbf{U}^\top(x,\lambda_j ) U_0(x,\lambda ))|_{\lambda=\lambda_j} p(x)d x.\end{eqnarray}
{This setup corresponds to the full upper triangular $K\times L$ \eqref{eq:ext.matrix} . However, as already mentioned, removal of any element of the data matrix  with the column number larger than K does not affect the computation of $\mathbf{U}$,  so we can use any sparsity patterns in the right  upper blocks of $(F_0-F)$ and $\frac{d}{d\lambda}(F_0-F)$.  Note that these patterns can also vary for different matching frequencies. This would result in the reduction of the range of the forward linear operator without changing the remaining columns.}

\section{Numerical experiments}
\subsection{One dimensional SIMO problem}
\begin{figure}[htb]
\centering
\includegraphics[scale=0.4]{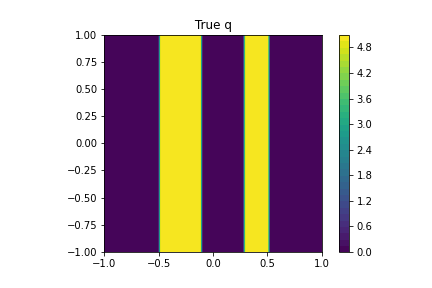}\includegraphics[scale=0.4]{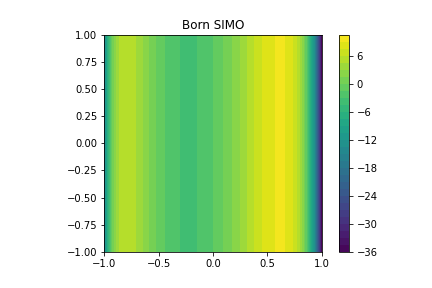}
\includegraphics[scale=0.4]{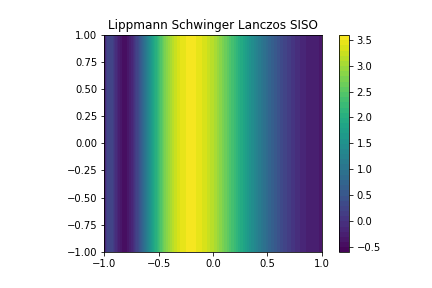}\includegraphics[scale=0.4]{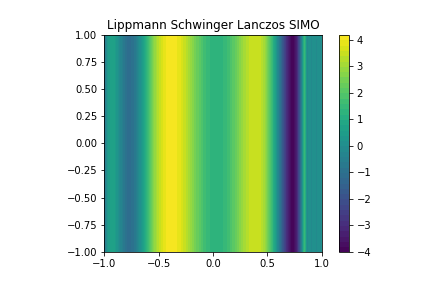} 
\caption{True one dimensional medium (top left) and its reconstructions using SIMO data and the Born approximation (top right), Lippmann Schwinger Lanczos using SISO data only (bottom left),  and Lippmann Schwinger Lanczos using SIMO data  (bottom right)}. 
\label{fig:num3}
\end{figure}

We first demonstrate numerical results here for a one dimensional reconstruction problem as a proof of concept. The setup is exactly as in Section 2, although computations are done on a square, using finite elements and $N=100\times 100$ elements to compute the synthetic data. The LSL system is solved here using the background solutions as a basis for the unknown $p$. We use $m=6$ positive spectral values $\lambda= 1,2, 14, 50, 128, 262.2672$.  Note in Figure \ref{fig:num3} that the SISO data, which is used to create the internal solutions, reconstructs the bar closer to the source quite well. Using the same internal solutions, we add data from the receiver on the right hand side, and we see that the second bar is captured.  The Born approximation does not work well for this high contrast example. 

\subsection{Two dimensional examples}

\begin{figure}[htb]
	\centering
	\begin{tabular}{cc}
		\includegraphics[width=0.45\textwidth]{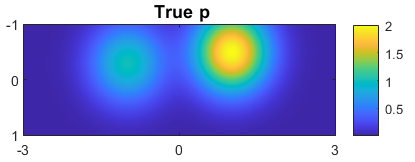} &
		\includegraphics[width=0.45\textwidth]{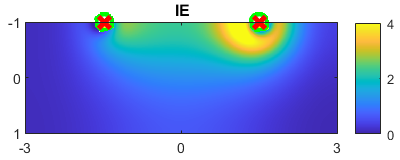}\\
		\includegraphics[width=0.45\textwidth]{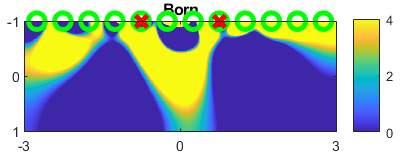} &
		\includegraphics[width=0.45\textwidth]{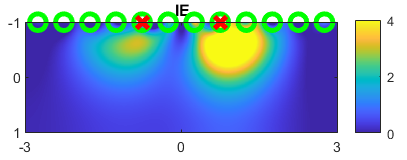}\\
	\end{tabular}
	\caption{{Experiment 1: True medium (top left) and its reconstructions using the LSL IE for square data $K=L=2$ (top right), Born linearization for rectangular data $K=2,~L=12$ (bottom left) and extended LSL IE for the same rectangular data(bottom right)}}
	\label{fig:num1}
\end{figure}

\begin{figure}[htb]
	\centering
	\hskip -.7in
	\includegraphics[width=0.45\textwidth]{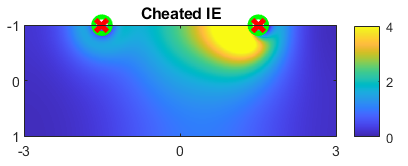}
	\includegraphics[width=0.45\textwidth]{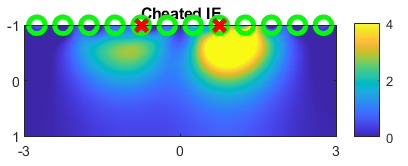}
	\caption{Experiment 1: {LSL IE  image with exact internal solution  ${U}(x,\lambda))$ (Cheated IE)} for square data $K=L=2$ (left) and for rectangular data $K=2,~L=12$ (right)}
	\label{fig:num2}
\end{figure}

Here we show our numerical reconstructions of a 2D two-bump media. We assumed that the data is noiseless. The construction of the data-driven ROM may become unstable in the presence of noise, however, one can regularize it using the SVD-based algorithm of \cite{borcea2019robust}. For simplicity, we omitted that part in this paper. We chose frequencies that provide enough sensitivity to recover the unknown scatterers,  but so that the data is still non-redundant. In our first experiment, we considered the medium shown on the top left in Fig.\ref{fig:num1}. The data acquisition setup mimics one from surface geophysical exploration, that is, we consider $K$ sources $L$ receivers on the top boundary of the domain. We use $m=5$ positive spectral values $\lambda=1,2,14,50,128$. The forward problem was discretized on a regular triangular grid using a finite-element method with $N=451\times 151$ elements. Equation (\ref{eq:intdMIMO2}) was approximated using nodal quadrature on a $901\times 301$ grid. The obtained linear system is ill-conditioned, and we solved it via projecting onto its dominant eigenvectors. On the top right of Fig.\ref{fig:num1} we plotted the reconstruction using the LSL IE for the square data matrix when $K=L=2$. Red crosses and green circles on plots show source(s) and receiver(s) locations, respectively. Then, we expanded the data matrix to the rectangular one by adding 10 receivers. For this scenario, the reconstruction using Born linearization, i.e. when ${U}(x,\lambda))$ in (\ref{eq:intdMIMO}) is replaced by with ${U}_0(x,\lambda))$), is shown on bottom left of Fig.\ref{fig:num1}. For the same data, the image produced by the extended LSL IE is plotted on the bottom right of Fig.\ref{fig:num1}. As one can observe, it  improves the image compared to both Born and the square LSL IE. In fact, in Fig.\ref{fig:num2} we plotted the solutions of (\ref{eq:intdMIMO}) assuming the exact internal solution ${U}(x,\lambda))$ is available, {which we call "Cheated IE". This represents the best possible result if one were to apply an iterative (aka distorted)  Lippmann-Schwinger algorithm which converged. (Iterative Lippman-Schwinger would require multiple solutions of forward problems, which can be costly even if such an algorithm converges.)  The cheated reconstructions are identical to the corresponding results for the LSL IE in Fig.\ref{fig:num1}, and, obviously, in both cases improvements are obtained compared to the rectangular data just because of the additional receivers.}

\begin{figure}[htb]
	\centering
	\hskip -.7in
	\includegraphics[width=0.9\textwidth]{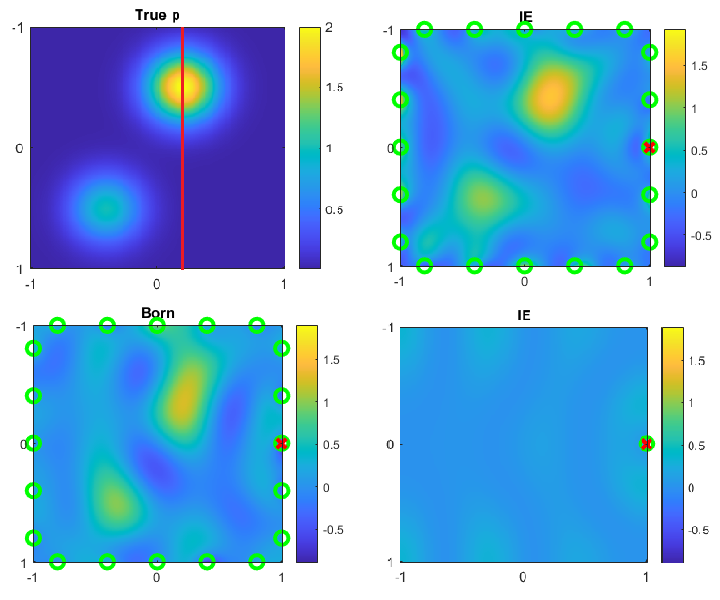}
	\caption{{True p with its LSL SIMO, Born SIMO, and LSL SISO images. The SIMO reconstructions are qualitatively correct, while the SISO reconstruction is not. The LSL SIMO (IE) is a quantitatively better approximation to the true $p$ than Born. }}
	\label{fig:num4}
\end{figure}

\begin{figure}[htb]
	\centering
	\hskip -.7in
	\includegraphics[width=0.9\textwidth]{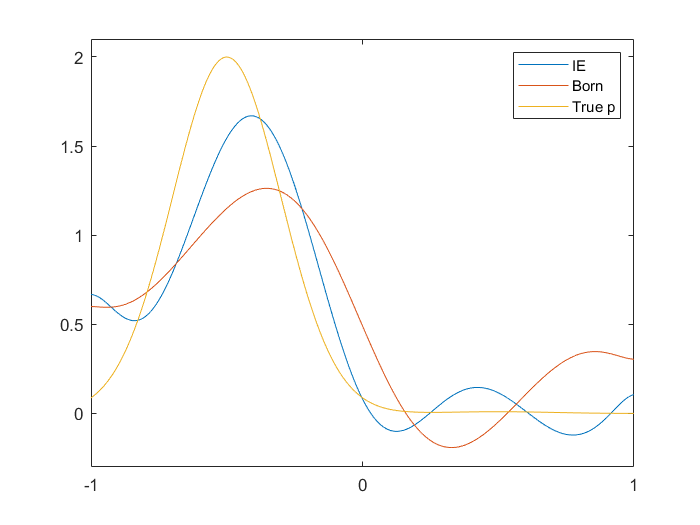}
	\caption{ Slices of true $p$, the SIMO Born and the  SIMO (LSL) IE.  The LSL IE reconstruction  is obviously closer to the true model.}
	\label{fig:slice}
\end{figure}

In the second experiment we considered a SIMO medical imaging setup with $K=1$ source and $L=20$ receivers located on all four sides of square domain (see top left in Fig.\ref{fig:num4}) and compared it to SISO scenario $K=L=1$. We discretized the problem using finite elements on a uniform triangular grid with $N=151\times 151$ elements. Equation \ref{eq:intdMIMO2} was approximated using quadrature on $451\times 451$ grid. In this example we considered $m=6$ spectral values $\lambda=-24,1,2,14,50,128$, that is, we added a negative one to the set from the previous example.  In Fig.\ref{fig:num4} we plotted the images obtained using the {LSL IE} for SIMO (top right) and SISO (bottom right) scenarios as well as the Born reconstruction for SIMO (bottom left). In this case the SISO LSL image is qualitatively wrong. Both SIMO LSL IE and born images are qualitatively correct, however the LSL IE image is quantitatively more accurate. {To strengthen this argument we compared a slice  (the red  line in Fig~\ref{fig:num4} top left ) of the true $p$ with its Born and LSL IE SIMO images in Fig.~\ref{fig:slice}.}


\section{Conclusion and discussion}
{Using the previously developed  Lippmann-Schwinger-Lanczos (LSL) algorithm, we are able to extend successfully the model reduction approach  to  arrays with nonreciprocal subsets, including SIMO data sets.  The main novelty of this algorithm is the exploitation of  the product structure of the  Lippmann-Schwinger equation. It allows us to  plug in an approximate internal solution computed using only the  data appropriate  for the ROM construction  (symmetric) subset of the data, and then using all of the given data to solve the enlarged linear  LSL  system.  Even though all of the derivations here are made for the second order Schrodinger equation in the Laplace domain, and assuming exact data, we plan to extended them to first order systems and the time domain, as well as regularized noisy ROM formulations of \cite{borcea2019robust,borcea2020reduced}. Future directions also include using iterative methods for very sparse data sets. We have already obtained promising preliminary results for  the monostatic (SAR) problem, which is an important case of such data sets.

\thanks{{\bf Acknowledgements} 
We are grateful Liliana Borcea, Alex Mamonov and J\"orn Zimmerling for productive discussions that inspired this research. V. Druskin was partially supported by AFOSR grant FA 955020-1-0079 {and NSF grant  DMS-2110773}. S. Moskow was partially supported by NSF grants DMS-1715425 and DMS-2008441.  

\bibliography{biblio,biblio6,graphbib6,galerkincitations}
\end{document}